\documentclass[12pt]{amsart}
\usepackage{amssymb,amsmath,amsthm,array}
\newtheorem{theorem}{Theorem}
\newtheorem{lemma}[theorem]{Lemma}
\newtheorem{corollary}[theorem]{Corollary}

\newtheorem{proposition}[theorem]{Proposition}
\newtheorem{definition}[theorem]{Definition}

\theoremstyle{remark}

\newtheorem*{remark}{Remark}

\numberwithin{theorem}{section} \numberwithin{equation}{section}

\newcommand{\HH}{\mathcal{H}}

\newcommand{\R}{\mathbb{R}}
\newcommand{\M}{\mathbb{M}}
\newcommand{\E}{\mathbb{E}}
\newcommand{\C}{\mathbb{C}}
\newcommand{\T}{\mathbb{T}}            

\newcommand{\Q}{\mathbb{Q}}

\newcommand{\Z}{\mathbb{Z}}
\newcommand{\N}{\mathbb{N}}

\def\H{\mathbb{H}}

\newcommand{\kzxz}[4]{\left( \begin{smallmatrix} #1 & #2 \\ #3 & #4\end{smallmatrix}\right) }
\newcommand{\kkzxz}[4]{\begin{smallmatrix} #1 & #2 \\ #3 & #4\end{smallmatrix} }

\newcommand{\nc}{\newcommand}
\nc{\cd}{{\mathcal D}}
\nc{\tildetau}{\tilde{\tau}}
\nc{\gauss}[1] {\lfloor #1 \rfloor} 
\nc{\gaussgauss}[2] {\lfloor \frac{#1}{#2} \rfloor}
\nc{\biggaussgauss}[2] {\bigg{\lfloor} \frac{#1}{#2} \bigg{\rfloor}}

\newcommand{\diff}[2]
{\cd_{#1,\underline{#2}}}
\nc{\diffknu}{\diff {k}{2 \nu}}
\nc{\dimens}{\mathsf{dim}}

\begin{document}
\title[A Strong Symmetry Property of Eisenstein series]{A Strong Symmetry Property of Eisenstein series}

\author{Bernhard Heim}
\address{Max-Planck Institut  f\"ur Mathematik, Vivaitsgasse 7, 53111 Bonn, Germany}
\email{heim@mpim-bonn.mpg.de}


\maketitle
\section{Introduction and Statement of results} 

Eisenstein series play a critical role in number theory. For two
hundred years they have been an essential
tool in the analysis of automorphic L-functions and in studying properties of
quadratic forms in one and several variables. 
The construction is clear and straightforward, while their
properties are sometimes very surprising. The arithmetic of their
Fourier coefficients, and their analytic properties are still not
completely understood. There are many connections with the Riemann
hypothesis and other famous unsolved problems in  
number theory.

Eisenstein series are named after Ferdinand Gotthold Eisenstein (1823 - 1852).
Let $k$ be an even integer larger than 2 and let $\tau$ be in the upper complex half-space.
One of the simplest Eisenstein series is defined by
\begin{equation}
E_k(\tau) := \frac{1}{2} \sum_{m,n \in \Z, \, (m,n)=1} \left( m \tau + n\right)^{-k}.
\end{equation}
It has the transformation property 
\begin{equation}
E_k\left( \frac{a \tau + b}{c \tau + d} \right) = \left( c \tau + d \right)^k E_k(\tau)
\end{equation}
for $ \left(\begin{smallmatrix}
a & b \\ c & d
\end{smallmatrix}\right) \in \text{SL}_2(\Z) $. It has a Fourier 
expansion with rational Fourier coefficients with bounded denominators,
involving divisor functions and Bernoulli numbers, and is
connected with special values of the Riemann zeta function. 

To understand special values of more general types of L-functions, this
simple version of Eisenstein series has been extended in many directions.
Siegel and Klingen studied Eisenstein series attached to the 
symplectic group, in order to study, for example, quadratic forms and the
structure of Siegel modular forms in several variables. 
This culminated in the Siegel-Weil formula \cite{We65} and the structure theorem. 
Later Klingen introduced the Eisenstein
series now called {\it Klingen type} \cite{Kl90}. 
In another direction, Maass, Roelcke, and Selberg
\cite{Se56} studied real analytic  
Eisenstein series in the context of differential operators and spectral theory.
Langlands \cite{La76} succeeded in showing remarkable general analytic
properties, i.e.,  
meromorphic continuation to the whole complex plane and functional equation,
for a wide range of reductive groups. This has applications in the
Rankin-Selberg and the Shahidi methods
to study analytic and arithmetic properties of automorphic
L-functions. The arithmetic properties of the Fourier  
coefficients play a fundamental role in the study of the
arithmetic of the special values. 
Garrett's integral representation of the triple L-function \cite{Ga87}
was an unexpected example of a different sort.

Yet another direction appears in the brilliant work of the late
H. Maass \cite{Ma79}, who found a new relation satisfied by the
Fourier coefficients of holomorphic Eisenstein series of Siegel type
of degree $2$. Automorphic functions with this property he called the
{\it Spezialschar}. 
His beautiful work on this subject made it possible to understand and prove 
the main part of the Saito-Kurokawa conjecture \cite{Za80}. 
Recently Skinner used results of Shimura on delicate properties
of Eisenstein series to attack the Iwasawa conjecture (see also \cite{Br07}). This brief review 
suggests that that new features of Eisenstein series should be fruitful.

In this paper we present a new method to study Fourier coefficients of 
holomorphic and non-holomorphic Eisenstein series simultaneously. 
This leads to a fundamental identity we state now.
We mainly focus on the real analytic Eisenstein series on Siegel upper
half-space $\H_2$ of degree $2$ to make our method clear and to not
burden the discussion with other technical considerations.

Let $E_k^{(2)}(Z,s)$ be the real analytic Eisenstein series 
of weight $k$ and $Z \in \H_2$ with respect to the Siegel modular group $\text{Sp}_2(\Z)$ 
and $s \in \C$ with $2 \, \text{Re}(s) + k >3$. For details we refer to section \ref{eisen}.
This function is not holomorphic as a function of $Z$ on $\H_2$, 
but does satisfy the transformation rule of a modular form. Since it
is periodic with respect to the real part $X$ of $Z$ it has a Fourier
expansion:   
\begin{equation}
E_k^{(2)}(X + i Y,s)=
\sum_{ N } A(N,Y;s) \, e^{2 \pi i \,\text{tr}\,(NX)}.
\end{equation}
where $Y$ is the imaginary part of $Z$ and $N$ is summed over 
half-integral matrices $N = \left( \kkzxz {n}{r/2}{r/2}{m}  \right)$.  

Then the following identity holds between the Fourier
coefficients $A(N,Y;s)$.  
Let $G[H] := H^t G H$ for appropriate matrices $G$ and $H$.
We have for all prime numbers $p$ and for all half-integral $N =\left( \kkzxz {n}{r/2}{r/2}{m}  \right)$ the formula
%
\begin{equation}
p^{k-1}
A  
\left(
\left(\begin{array}{cc}
\frac{n}{p}& \frac{r}{2p} \\ \frac{r}{2p}& m \end{array}\right), p\, Y;s
\right) 
-A \left(
\left(\begin{array}{cc}
n & \frac{r}{2} \\ \frac{r}{2}& p\, m \end{array}\right),
Y
;s
\right)
\end{equation}
\begin{eqnarray*}
& = & p^{k-1}
A  
\left(
\left(\begin{array}{cc}
n& \frac{r}{2p} \\ \frac{r}{2p}& \frac{m}{p} \end{array}\right), 
Y 
\left[
\left(\begin{array}{cc}
1& 0 \\ 0& p \end{array}\right)
\right]
;s
\right) \\
& & \,\,- \,\,\,\,
A \left(
\left(\begin{array}{cc}
pn & \frac{r}{2} \\ \frac{r}{2}& m \end{array}\right),
Y\left[
\left(\begin{array}{cc}
p^{-\frac{1}{2}}& 0 \\ 0& p^{\frac{1}{2}} \end{array}\right)\right];s
\right).
\end{eqnarray*}
Here we put $A(N,Y,s)=0$ if $N$ is not half-integral.

The nature of the Fourier coefficients $A(N,Y;s)$ is complicated,
involving special values of Dirichlet L-series (Siegel series) and
Bessel functions of higher order.  
One has to distinguish the various cases of the rank of $N$. 
Nevertheless, our method works without any explicit knowledge of 
these formulas, and is completely explicit. Moreover it also works in
the case of Hecke summation. 

This paper is organized in the following way. In $\S2$ we recall some
basic aspects of Shimura's approach to the theory of Hecke. 
This will be used to define a new kind of operators, which do not act
on the space of modular forms, but nevertheless 
inherit interesting properties. In $\S3$ we prove a decomposition of
the real-analytic Eisenstein series, essentially based on
consideration of two subseries $A_k(Z,s)$ and $B_k(Z,s)$, concerning
which we prove several properties. In $\S4$ we present the main result
of this paper, namely, we show that real-analytic Eisenstein series
satisfy the strong symmetry property 
\begin{equation}
E_k^{(2)}(Z,s) \vert \bowtie T = 0
\end{equation} for all Hecke operators $T$, 
which will be explained in that paragraph in detail.
We also give an example of a family of modular forms which do not have this property.
Finally, we give applications, for example, the fundamental identity among
the Fourier coefficients. 
\section{Hecke Theory 
\'{a}
la Shimura}
For $k \in \N$ be even let $M_k$ be the space of elliptic modular forms of weight $k$ 
with respect to the full modular group $\Gamma = \text{SL}_2(\Z)$.
Let $f \in M_k$. Hecke introduced the operators
$T_n, n \in \N$ given by
\begin{equation}
T_n(f)(\tau) := n^{k-1} \sum_{d \vert n} d^{-k} 
\sum_{b=0}^{d-1} f \left( \frac{n \tau + bd}{d^2}\right),
\end{equation}
which map modular forms to modular forms. These operators commute with
each other. They are multiplicative 
and self-adjoint with respect to the Petersson scalar product on the space of cusp forms. The vector space $M_k$ has a 
basis of simultaneous eigenforms. The eigenvalues $\lambda_n(f)$ are totally real integers and are proportional 
to the $n$-th Fourier coefficients of the eigenform. Shimura
\cite{Sh71} studied systematically  
the underlying Hecke algebra. The realization of this Hecke algebra on
the space of modular forms gives then the Hecke operators above. 
\\
\\
We start with some basic constructions \cite{Sh71}. Let $(R,S)$ be a 
Hecke pair, meaning that $R$ is a subgroup of the group $S$
and for each $s \in S$ the coset space $R \backslash RsR$ is finite. 
For $P$ be a principal ideal domain, $R$ acts on the right on the 
$P$-module $L_P(R,S)$
of formal finite sums $X = \sum_{j} a_j R s_j$ with $a_j \in P, s_j \in S$. 
The subset $H_P(R,S)$ of elements invariant under this 
action forms a ring with the multiplication
\begin{equation}
\left( \sum_i a_i R g_i \right) \circ 
\left( \sum_j b_j R h_j \right)  := \sum_{i,j} a_i b_j R g_i h_j. 
\end{equation}
This ring is called the associated Hecke ring or algebra.
It is convenient to identify the left coset decomposition of the double cosets 
$RsR = \bigsqcup_j R s_j$ with $\sum_j Rs_j \in H_P(R,S)$
which form a basis of the $P$-module $H_P(R,S)$. 
Hence double cosets are identified with a full system of representatives of
the $R$-left coset decomposition of the double coset.

Now we apply this construction to our situation. For $l \in \N$ put
\begin{eqnarray}
M(l)&:= &
\bigsqcup_{d\vert l, \, \, d\vert \frac{l}{d}} 
\Gamma \kzxz{d}{0}{0}{l/d} \Gamma. 
\end{eqnarray}
Then we set $M_{\infty} := \bigsqcup_{l=1, n = 1}^{\infty} 
\left(\begin{smallmatrix}
n^{-1} & 0 \\ 0 & n^{-1}
\end{smallmatrix}\right) 
M(l)$. The following property is well-known.
					\begin{lemma}
					We have  that
					$\left( \Gamma, M_{\infty}\right) $ is a Hecke pair.
					\end{lemma}
Let ${\HH}$ be the corresponding Hecke algebra of the Hecke pair
$\left( \Gamma, M_{\infty}\right) $ over $\Q$. 
Then we have the Hecke pair
$$
\Big( \Gamma \, , \,\cup_{l \in \Z} M(p^l)            \Big) 
$$
for all prime $p$ with corresponding Hecke algebra $\HH_p$.
%
By the elementary divisor theorem
\begin{equation}
\HH = \otimes_p \HH_p.
\end{equation}
Let $\T_l= \Gamma \backslash M(l)$. Then
the Hecke algebra $\HH_p$ 
is generated by the $\T_p$, the special double cosets
$
\Gamma \kzxz{1}{0}{0}{p} \Gamma
$, and $\Gamma
\left(\begin{smallmatrix}
p^{-1} & 0 \\ 0 & p^{-1}
\end{smallmatrix}\right) 
\Gamma$. 
Here
\begin{equation}
\Gamma\kzxz{1}{0}{0}{p} \Gamma = 
\Gamma \kzxz{p}{0}{0}{1}  + \sum_{a=0}^{p-1}\Gamma \kzxz{1}{a}{0}{p}. 
\end{equation}
Let $\text{GL}_2^{+}(\R)$ the set of $\R$-valued $2\times 2$ matrices with positive determinant.  Let $M \in \text{GL}_2^{+}(\R)$. 
Define $\widetilde{M} := \text{det}(M)^{-\frac{1}{2}} \, M$.
\begin{definition}
The action of the Hecke algebra $\HH$ on 
$M_{k}$ is induced by double cosets. Let $g \in \text{ \rm GL}_2^{+}(\Q)$ and $f \in M_k$. Then 

\begin{equation}
f \vert_k [ \Gamma g \Gamma] := \sum_{ A \in \Gamma \backslash 
\Gamma g \Gamma
} f\vert_k \widetilde{A}.
\end{equation}
Here $\vert_k$ is the Petersson slash operator. In particular, the 
normalized Hecke operators are defined by
\begin{equation}
\T_n (f) := n^{\frac{k}{2}-1} 
\sum_{ A \in \Gamma \backslash M(n)} f\vert_k \widetilde{A}.
\end{equation}
\end{definition}
\remark
The Hecke operators $\T_n$ coincide with the classical Hecke operators 
$T_n$ on the space $M_k$. For $f$ be a primitive form, the eigenvalue
of $T_n$ is the $n$-th Fourier coefficient of $f$.
\\
\\
Shimura's approach to Hecke theory can be generalized to introduce new operators 
related to classical Hecke operators, and which coincide in certain
special situations. 

Let $ A = \kzxz {a}{b}{c}{d}$ and $B = \kzxz {e}{f}{g}{h}$. Then 
\begin{equation} A \times B := 
\kzxz
{\kkzxz {a}{0}{0}{e}} 
{\kkzxz {b}{0}{0}{f}} 
{\kkzxz {c}{0}{0}{g}} 
{\kkzxz {d}{0}{0}{h}}
\end{equation}
gives an embedding of $\text{SL}_2(\R) \times \text{SL}_2(\R)$ into the symplectic group $\text{Sp}_2(\R)$ of degree $2$.
Let $A \in \text{GL}_2(\R)$ with det($A$)$= l >1$. We put
\begin{equation}
\widetilde{A}^{\bullet} := \kzxz{l^{-1/2} a}{\,\,l^{-1/2} b}{l^{-1/2} c}{\,\,l^{-1/2} d} \times 
\kzxz{1}{0}{0}{1},
\end{equation}
and similarly define $\widetilde{A}_{\bullet}$.
Let $F: \H_2 \longrightarrow \C$ with $F\vert_k g^{\bullet} = F$ for all $g\in \Gamma$. Let
$A \in \text{GL}_2^{+}(\Q)$.
Define the Hecke operator
\begin{equation}
F\vert_k \widetilde{[ \Gamma A \Gamma]}^{\bullet} := 
\sum_{\gamma \in \Gamma \backslash \Gamma A \Gamma} F\vert_k \widetilde{\gamma}^{\bullet},
\end{equation}
and similarly
$F\vert_k \widetilde{[ \Gamma A \Gamma]}_{\bullet} $. 
For simplicity put $\vert_k T^{\bullet}$ and $\vert_k T_{\bullet}$ for $ T \in \HH$.
\section{Eisenstein series decompositon}\label{eisen}
In this section we state and prove a 
decomposition formula for $E_k^{(2)}(Z,s)$.
It is essentially constructed from two functions. 
The symplectic group $\text{Sp}_n(\R)$ acts on the Siegel 
upper half-space $\H_n$ of degree $n$ via 
$
\left(\begin{smallmatrix}
A & B \\ C & D
\end{smallmatrix}\right) 
(Z ):= \left( AZ + B \right) \left( CZ + D \right)^{-1}$.
We put $j\left(
\left(\begin{smallmatrix}
A & B \\ C & D
\end{smallmatrix}\right),Z 
\right):= \text{det}\left( CZ +D \right)$. 
Let $\Gamma_n:= \text{Sp}_n(\Z)$ be the Siegel modular group 
and let $\Gamma_{n,0}$ be the subgroup of all elements with $C=0$.
					\begin{definition}
					Let $k$ be an even integer and let $n \in \N$. Define the real 
					analytic Eisenstein series of weight $k$ and genus $n$ on
					$\H_n \times \mathcal{D}_k^n$, where $$\mathcal{D}_k^n := \{ s \in \C \vert \, 2 \,
					\text{ \rm Re}(s) + k > n +1\},$$ by
					\begin{equation}\label{infinite}
					E_k^{(n)} \left( Z,s\right) := \sum_{ g \in \Gamma_{n,0} 
					\backslash \Gamma_n} j(g,Z)^{-k} \, \delta\left( g(Z)\right)^{s}.
					\end{equation}
					Here $\delta\left( Z \right) := \text{ \rm det}\left( \text{ \rm Im}(Z)\right)$.
					\end{definition}
The infinite sum in (\ref{infinite}) converges absolutely and
uniformly on compacts
on the set $\H_n \times \mathcal{D}_k^n$. From Langlands' theory \cite{La76},
$E_k^{(n)}(Z,s)$ has a meromorphic continuation in $s$ to the whole
complex plane, and satisfies a functional equation.
In particular, let $k$ be an even positive integer, let $\xi(s) :=
\pi^{\frac{s}{2}} \Gamma(\frac{s}{2})  
\zeta(s)$ and $\Gamma_n(s) := \prod_{j=1}^{n} \Gamma(s - \frac{j-1}{2})$.
Here $\Gamma(s)$ is the Gamma function and $\zeta(s)$ the Riemann zeta
function. Then the function 
\begin{equation}
\E_k^{(n)}(Z,s) := \frac{\Gamma_n(s + \frac{k}{2})}{\Gamma_n(s)} \cdot\xi(2s) 
\prod_{i=1}^{[n/2]} \xi(4s-2i) \,\, E_k^{(n)}\left( Z, s-\frac{k}{2}\right)
\end{equation}
is invariant under $s \mapsto \frac{n+1}{2} -s$ and is an entire
function in $s$ (see \cite{Mi91}). 
Here $[x]$ is the largest integer smaller or equal to $x$. When $n=1$ the function
\begin{equation}
\E_k(\tau,s) = \Gamma\left(s + \frac{k}{2}\right) \, \zeta(2s) \, \pi^{-s} \, E_k\left(\tau, s - \frac{k}{2}\right)
\end{equation}
is entire and is invariant under $s \mapsto 1-s$. Moreover, for $n=2$ the function 
\begin{eqnarray*}
\E_k^{(2)}(Z,s) &= &\Gamma(s) \, \Gamma\left(s+\frac{k}{2}\right) \, \Gamma\left( s + \frac{k-1}{2}\right) \,2^{2s-2} \,
\pi^{-s - \frac{1}{2}} \\ 
& &\zeta(2s) \,\zeta(4s -2) \,E_k^{(2)}\left(Z, s -\frac{k}{2}\right)
\end{eqnarray*}
entire and invariant under $s \mapsto \frac{3}{2} -s$.
\\
\\

For a positive even integer $k$ with $k > n+1$ the function 
$E_k^{(n)}(Z) := E_k^{(n)}(Z,0)$ is the holomorphic Siegel Eisenstein series.
It has a Fourier expansion with rational coefficients. 
Moreover the denominators are bounded. In the real analytic case the
situation is somehow different. 
The Fourier coefficient depend on the imaginary part of $Z$ and
involve confluent hypergeometric functions. 
Moreover, one has to study Hecke summation if one is interested in the
case $k = n+1$ and $s=0$, for example. 
Let $k$ be an even integer. Then
$\mathcal{D}_k := \{ s \in 
\C\vert \, 2 \, Re (s) + k > 3\}$. It is well known that $E_2^{(2)}(Z,0)$ is finite. 
But we do not want to go into this topic further.
We parametrize $Z \in \H_2$ by
$\left(\begin{smallmatrix}
\tau & z \\ z & \tildetau
\end{smallmatrix}\right)$ and define
$
\varphi_k (Z) := \tau + 2z + \tildetau,
$.
For simplicity, put 
$\chi_{k,s}(g,Z) := j(g,Z)^{-k} 
\vert j(g,Z)\vert^{-2s}$ and 
$\Phi_{k,s} := \varphi_k(Z)^{-k} 
\vert \varphi_k(Z)\vert^{-2s}$ for $g \in \text{Sp}_2(\R)$.
Also let $\Gamma_{\infty} = \Gamma_{1,0}$ and $\H = \H_1$.
Let $\vert_k$ be the Petersson slash operator. We drop the
symbol for the weight $k$ if it is clear from the context.
\begin{definition}
For $k \in \Z $ be even we define two $\C$-valued functions $A_k$ (resp. $B_k$) 
on $\H_2 \times \mathcal{D}_k$
by
\begin{eqnarray*}
\left( Z,s \right)   &\mapsto &
\delta(Z)^s \sum_{ g, h \in \Gamma_{\infty}\backslash \Gamma}\limits
\chi_{k,s}\left(g^{\bullet}h_{\bullet},Z \right) \quad
\text{ and } \\
\left( Z,s \right)  & \mapsto &
\delta(Z)^s \sum_{g \in \Gamma}\limits \Phi_{k,s}\left( 
g_{\bullet}(Z)\right) \chi_{k,s} \left( g_{\bullet},Z\right).
\end{eqnarray*}
\end{definition}
These functions turn out to be subseries of the real analytic
Eisenstein series of degree two, with similiar convergence properties.
					\begin{theorem}
					Let $k$ be an even integer. Let $Z \in \H_2$ and 
					$s \in \mathcal{D}_k$. Then
					\begin{equation}\label{decomposition}
					E_k^{(2)} \left( Z,s\right) = A_k(Z,s) + 
					\sum_{m=1}^{\infty} B_k\vert
					\left(
					\Gamma 
					\left( 
					\kkzxz {m}{0}{0}{m^{-1}} \right)\Gamma \right)
					^{\bullet}
					\left( Z,s\right)\, \, m^{-2s -k}.
					\end{equation}
					\end{theorem}
\begin{proof}
From Garrett \cite{Ga84}, \cite{Ga87} we know how to study coset systems of the type
$$
\Gamma_{2n,0} \backslash \Gamma_{2n}/ \Gamma_n \times \Gamma_n
$$
in the context of the doubling method. Similarly, we
obtain a useful $\Gamma_{2,0}$-left coset decomposition of $\Gamma_2$
given by $R_0 \bigsqcup R_1$ with  
\begin{equation}
R_0 = \Gamma_{\infty} \backslash \, \Gamma
\times
\Gamma_{\infty} \backslash \, \Gamma \text{    and    }
R_1 = \bigsqcup_{m=1}^{\infty} g_m \, \Big( \Gamma \times \Gamma(m)\backslash \Gamma \Big).
\end{equation}
Here $\Gamma(m) := \{ g \in \Gamma \vert \, \kzxz{0}{1/m}{m}{0} g \kzxz{0}{1/m}{m}{0} \in \Gamma\}$ and
\begin{equation}
g_m := \left(\begin{array}{cccc}
 0& 0& -1&0 \\
0&1 & 0& 0\\
1&m & 0&0 \\
0&0 & -m& 1
\end{array}\right).
\end{equation}
The subseries related to the representatives $g_m \left(\Gamma
  \times \Gamma(m)\backslash \Gamma \right)$ is 
\begin{equation}\label{formula1}
\delta(Z)^s \sum_{ g \in \Gamma,\, h \in \Gamma(m)\backslash \Gamma }
\chi_{k,s} \left( g_m \, (g\times h), Z\right).
\end{equation}
Let ${\M}_{m} $ be the diagonal $4\times 4$ matrix 
with $(1,m,1,m^{-1})$ on the diagonal. Then
$j(g_m,Z) = j\left(g_1, \M_m(Z)\right)$. 
Hence we obtain, for (\ref{formula1}):
\begin{equation*}
\delta(Z)^s \sum_{ g \in \Gamma, g \in \Gamma(m)\backslash \Gamma }
\Phi_{k,s} \left( \M_m \, (g \times h)(Z)\right) \chi_{k,s} \left(g \times h,Z\right).
\end{equation*}
Let $\#$ be
the automorphism of $\text{SL}_2(\R)$ given by
$\left(\begin{smallmatrix}
a & b \\ c & d
\end{smallmatrix}\right)^{\#} := \left(\begin{smallmatrix}
d & b \\ c & a
\end{smallmatrix}\right)$ of $\text{SL}_2(\R)$. Then we can prove in
a straightforward manner the symmetric relation
\begin{equation}
\Phi_{k,s} \left( g^{\bullet}(Z) \right) \chi_{k,s}(g^{\bullet},Z) =
\Phi_{k,s} \left( g_{\bullet}^{\#}(Z) \right) \chi_{k,s}(g_{\bullet}^{\#},Z) .
\end{equation}
By the elementary divisor theorem we obtain for our subseries the expression
\begin{equation}
\delta(Z)^s \, m^{(k+2s)} \sum_{ \gamma \in 
\Gamma 
\left( 
\kkzxz {m}{0}{0}{m^{-1}} \right)\Gamma} 
\Phi_{k,s} \left( 1_2 \times \gamma)(Z)\right) \, \chi_{k,s} \left(1_2 \times \gamma,Z\right).
\end{equation}
Now we can apply again the symmetry relation and obtain the formula
(\ref{decomposition}) in our theorem. 
\end{proof}
\begin{corollary}
Let $k$ be an even integer. Let $Z \in \H_2$ and 
					$s \in \mathcal{D}_k$. Then
					\begin{equation}
					E_k^{(2)} \left( Z,s\right) = A_k(Z,s) + 
					\sum_{m=1}^{\infty} \Big(B_k\vert
					\left(
					\Gamma 
					\left( 
					\kkzxz {m}{0}{0}{m^{-1}} \right)\Gamma \right)
					_{\bullet}
					\left( Z,s\right) \Big) \, m^{-2s -k}.
					\end{equation}
\end{corollary}
Let $F$ be a complex-valued function on $\H_2$. Let $k \in \N_0$ be even. 
Then we say that $F$ is $\Gamma$-modular of weight $k$ if
$F\vert_k \gamma^{\bullet} = F_k \vert 
\gamma_{\bullet} = F$ for all $\gamma \in \Gamma$.
\begin{corollary}
The functions $A_k(Z,s)$ and $B_k(Z,s)$ are $\Gamma$-modular.
\end{corollary}
\section{strong symmetry of Eisenstein series}
Let $F$ be a complex valued $C^{\infty}$ function on 
the Siegel upper half-space of 
degree $2$ with the transformation property of a modular 
form of even weight $k$ with respect to $\text{Sp}_2(\Z)$.
Let $f(\tau,\tildetau) := F \left( \kkzxz {\tau}{0}{0}{\tildetau} \right)$. 
Then we have the symmetry
\begin{equation}
f(\tau, \tildetau) = f(\tildetau, \tau),
\end{equation}
since $F\vert U =F$ with
$$
U:= \left(\begin{array}{cccc}
 0& 1& 0&0 \\
1&0 & 0& 0\\
0&0 & 0&1 \\
0&0 & 1& 0
\end{array}\right).
$$
It is worth noting that this does {\it not} imply 
that if we apply Hecke operators $T \in \HH$ on $f$ by 
fixing one of the variables that such a symmetry still holds.
Let for example $F$ be the holomorphic Klingen Eisenstein series of
degree $2$ and weight $12$ attached to the 
Ramanujan $\Delta$-function. Then it can be shown that
\begin{equation}
f(\tau, \tildetau) = E_{12}(\tau) \Delta(\tildetau) + E_{12}(\tildetau) \Delta(\tau) + \alpha \Delta(\tau) \Delta(\tildetau),
\end{equation}
with $\alpha \in \C$. Since infinitely many Hecke eigenvalues of the
Eisenstein series and the $\Delta$ function are different, it is
obvious that
\begin{equation}
f\vert \widetilde{T_p}^{\bullet} - f\vert \widetilde{T_p}_{\bullet} \neq 0
\end{equation}
for (at least) one prime number $p$.

The real analytic Eisenstein series $E_k^{(2)}(Z,s)$ of degree two 
has an important symmetry which had not been discovered
before. Let $T$ be an element of the Hecke algebra $\HH$. 
We will show in this section that, if we apply $T$ as an operator on the 
Eisenstein series to the two embeddings 
$T^{\bullet}$ and $T_{\bullet}$ we get the same new function, i.e.,
\begin{equation}
\left( E_k^{(2)}\vert \widetilde{T}^{\bullet}\right) \left(Z,s \right) = 
\left( E_k^{(2)}\vert \widetilde{T}_{\bullet}\right) \left(Z,s \right) .
\end{equation}
From the viewpoint of physics this can been seen as a scattering
experiment with an object $X$, in which we hit the object from outside
with $T_p$ for different prime numbers and look at the reaction. 
For example, if we knew in advance that the object were a holomorphic
Eisenstein series, then we could conclude that it is of Siegel type.

Actually we show that the subseries $A_k(Z,s)$ and
\begin{equation}
B_k^m(Z,s) := B_k \vert
					\left(
					\Gamma 
					\left( 
					\kkzxz {m}{0}{0}{m^{-1}} \right)\Gamma \right)
					_{\bullet}
					\left( Z,s\right)
\end{equation}
already have the strong symmetry property. Further, the function $A_k(Z,s)$ turns out
to be an eigenfunction.

					\begin{proposition}
					Let $k$ be an even integer and $s \in {\mathcal{D}}_k$. 
					For $T \in \HH$ we have
					\begin{equation}
					\left( A_k\vert \widetilde{T}^{\bullet} \right)(Z,s) =
					\left( A_k\vert \widetilde{T}_{\bullet} \right) (Z,s) 
					=  \lambda(T) A_k(T,s),
					\end{equation}
					with $\lambda(T) \in\C$.
					\end{proposition}
\begin{proof}
We have that
\begin{equation*}
A_k(Z,s) = \sum_{g,h \in \Gamma_{\infty} \backslash \Gamma} 
j( g^{\bullet}h_{\bullet},Z)^{-k} \delta\left( g^{\bullet}h_{\bullet}(Z))\right)^s.
\end{equation*}
At this point we note that $g^{\bullet}h_{\bullet} = h_{\bullet} g^{\bullet}$ 
and $j(g^{\bullet} h_{\bullet},Z) = j(g^{\bullet},h_{\bullet}(Z)) j(h_{\bullet},Z)$.
Since the series convergences absolutely and uniformly on compacts in 
$\H_2 \times \mathcal{D}_k$ we can interchange summation to obtain
\begin{eqnarray*}
A_k(Z,s) & = & \sum_{ h \in \Gamma_{\infty} \backslash \Gamma} j(h_{\bullet},Z)^{-k} 
\sum_{ g \in \Gamma_{\infty} \backslash \Gamma} j(g^{\bullet},h_{\bullet}(Z))^{-k}
\delta( g^{\bullet}(h_{\bullet}(Z))^{s} \\
& = & \sum_{ h \in \Gamma_{\infty} \backslash \Gamma} E_k \left( (h_{\bullet}(Z))^{*},s \right) j(h,Z_{*})^{-k}.
\end{eqnarray*}
Let $ Z = \kzxz {\tau}{z}{z}{\tildetau}$. Here $Z^{*} := \tau$ and
$Z_{*} := \tildetau$. By the same procedure we obtain
\begin{equation*}
A_k(Z,s) = \sum_{ g \in \Gamma_{\infty} \backslash \Gamma} E_k \left( (g^{\bullet}(Z))_{*},s \right) j(g,Z^{*})^{-k}.
\end{equation*}
Now let $ T \in \HH$ and $T = \sum_j a_j \Gamma t_j$. Then we have
\begin{equation*}
\left(A_k\vert \widetilde{T}^{\bullet}\right) (Z,s) = \sum_j a_j\sum_{ h \in \Gamma_{\infty} \backslash \Gamma}
E_k \left( (h_{\bullet}\widetilde{t_j}^{\bullet}(Z))^{*},s \right) j(h,\widetilde{t_j}(Z)_{*})^{-k} 
j(\widetilde{t_j}^{\bullet}, \tildetau)^{-k}.
\end{equation*}
Hence,
\begin{eqnarray}
\left( A_k \vert \widetilde{T}^{\bullet} \right) (Z,s) & = & 
\sum_j a_j\sum_{ h \in \Gamma_{\infty} \backslash \Gamma}
E_k \left( (\widetilde{t_j}^{\bullet}h_{\bullet}(Z))^{*},s \right)
j\left( \widetilde{t_j}^{\bullet}, h_{\bullet}(Z)^{*}\right)^{-k}
j \left( h_{\bullet},Z \right)^{-k} \nonumber\\
& = & 
\sum_{ h \in \Gamma_{\infty} \backslash \Gamma}
\left(E_k \vert \widetilde{T}^{\bullet}\right) 
\left( h_{\bullet}(Z)^{*},s\right) \, j\left(h_{\bullet},Z\right)^{-k}.
\end{eqnarray}
It is well known that $E_k(\tau,s)$ with $\tau \in \H$ is a Hecke eigenform. This leads to 
$\left( A_k \vert \widetilde{T}^{\bullet} \right) (Z,s) = \lambda(\widetilde{T}) \, A_k (Z,s)$. The same argument works for 
$\left( A_k \vert \widetilde{T}_{\bullet} \right) (Z,s)$ with the same eigenvalue. This proves the proposition. 
\end{proof}
\begin{proposition}
Let $k$ be an even integer. Let $m \in \N$ and let $T \in \HH$. Then we have             
\begin{equation}
\left( B_k^m\vert \widetilde{T}^{\bullet} \right)(Z,s) =\left( B_k^m\vert \widetilde{T}_{\bullet} \right) (Z,s)
\end{equation}
for all $\left( Z,s \right) \in \H_2 \times \mathcal{D}_k$.
\end{proposition}
\begin{proof}
Let $T = \sum_j a_j \Gamma g_j$ with $a_j \in \C$ and $g_j \in \text{Gl}_2^{+}(\Q)$.
Then we have
\begin{eqnarray*}
\left( B_k^m \vert \widetilde{T}^{\bullet} \right) (Z,s) & = & 
\sum_j a_j \,B_k \vert \left(
					\Gamma 
					\left( 
					\kkzxz {m}{0}{0}{m^{-1}} \right)\Gamma \right)_{\bullet} \widetilde{g_j}
					_{\bullet}(Z,s)\\
					& = & 
\sum_j a_j \,B_k \vert  \widetilde{g_j}^{\bullet}
					\left( \Gamma 
					\left( 
					\kkzxz {m}{0}{0}{m^{-1}} \right)\Gamma   \right)
					_{\bullet}(Z,s)
\end{eqnarray*}
since the Hecke algebra $\HH$ is commutative. Hence we can reduce our
calculations to the case $m=1$. Then we have for $\left(B_k  \vert
  \widetilde{T}^{\bullet} \right)(Z,s)$ the expression 
\begin{equation*}
\sum_j a_j \,\delta\left(\widetilde{g_j}^{\bullet}(Z)\right)^s \sum_{ g \in 
\Gamma} \Phi_{k,s} \left( (g_{\bullet} \widetilde{g_j}^{\bullet})(Z)\right) 
\chi_{k,s} \left(g_{\bullet} , \widetilde{g_j}^{\bullet}(Z)\right) \,\, j(g_j^{\bullet},Z)^{-k}.
\end{equation*}
To proceed further we use the cocycle property
\begin{equation*}
\chi_{k,s} \left( g_{\bullet} \widetilde{g_j}^{\bullet},Z \right) = 
\chi_{k,s} \left( g_{\bullet}, \widetilde{g_j}^{\bullet}(Z) \right) 
\chi_{k,s} \left( \widetilde{g_j}^{\bullet},Z \right)
\end{equation*}
and the transformation property $\delta \left( \widetilde{g_j}^{\bullet}(Z) \right)^s = 
\delta(Z)^s \vert j( \widetilde{g_j}^{\bullet},Z \vert^{-2s}$.
Hence $\left(B_k  \vert \widetilde{T}^{\bullet} \right)(Z,s)$ is equal to
\begin{equation*}
\sum_j a_j \delta(Z)^s \sum_{ g \in \Gamma} 
\Phi_{k,s} \left( (g_{\bullet} \widetilde{g_j}^{\bullet})(Z)\right) 
\chi_{k,s}\left( (g_{\bullet} \widetilde{g_j}^{\bullet}),Z\right) .
\end{equation*}
Now we apply the symmetry relation and note that $T$ is invariant with
respect to the automorphism $\#$. 
Then we obtain
\begin{equation*}
\sum_j \, \delta(Z)^s \sum_{g \in \Gamma} 
\Phi_{k,s} \left( ( \widetilde{g_j}_{\bullet}g_{\bullet})(Z)\right) 
\chi_{k,s}\left( (\widetilde{g_j}_{\bullet}g_{\bullet} ),Z\right) .
\end{equation*}
Finally we use the $\Gamma$-invariance property of $\Phi_{k,s}$ and $\chi_{k,s}$. This leads to
\begin{equation}
\sum_j \, \delta(Z)^s \sum_{g \in \Gamma} 
\Phi_{k,s} \left( ( g_{\bullet}\widetilde{g_j}_{\bullet})(Z)\right) 
\chi_{k,s}\left( (g_{\bullet}\widetilde{g_j}_{\bullet}),Z\right) .
\end{equation}
This gives the proposition.
\end{proof}
For $T \in\HH$ and even integer $k$ let $\vert_k \bowtie T$ be the 
operator $\vert_k \widetilde{T}^{\bullet} - \vert_k \widetilde{T}_{\bullet}$. If a $\Gamma$-modular function 
is annihilated by this operator, we say that is satisfies the strong symmetry property. 
This makes sense since this property turns out to classify
certain subspaces and gives a fundamental identity between Fourier coefficients.
Summarizing our results, we have
\begin{theorem}
Let $k$ be an even integer. Let $T$ be an element of the 
Hecke algebra $\HH$. Let $\left(Z,s\right) \in \H_2 \times \mathcal{D}_k$. Then we have
\begin{equation}\label{super}
E_k^{(2)}\vert \bowtie T (Z,s) = 0.
\end{equation}
\end{theorem}
\begin{corollary}
The strong symmetry (\ref{super}) of the Eisenstein 
series is also preserved under meromorphic continuation.
\end{corollary}
It would be interesting to study the implication of this property 
for the residues in relation with the Siegel-Weil formula.

\section{Applications of the strong symmetry property}
%
%
In \cite{He06} we have shown that a Siegel modular
form $F$ of degree $2$ 
with respect to the Siegel modular group $\text{Sp}_2(\Z)$ 
$F$ is a Saito-Kurokawa lift if and only if $F$ has the strong symmetry property.
Moreover, this can be used to study the non-vanishing of certain special values
predicted by the Gross-Prasad conjecture and in 
the context of the Maass-Spezialschar results recently proven by
Ichino.
Our proof in the holomorphic case was based on the interplay between Taylor coefficients and
certain differential operators. In this paper in the setting of real analytic Eisenstein series
the proof does not work. That was the reason why we gave a new one and which works
just because of the definition of an Eisenstein series via certain left cosets.
\begin{theorem}
Let $k$ be an even integer. Let $F: \H_2 : \longrightarrow \C$ be a 
$\C^{\infty}$-function which satisfies the transformation law
$
F \vert_k \gamma = F \text{  for all }\gamma \in \Gamma_2$.
Then we have
\begin{eqnarray}
F\vert \bowtie_{T}             \!\! &= &   \!\!          0 \text{ for all }           T \in \HH    \label{eins}     \\
& \Longleftrightarrow & \nonumber \\
F\vert \bowtie_{T_p}          \!\! &= &   \!\!          0 \text{ for all prime numbers } p \label{zwei}\\
& \Longleftrightarrow & \nonumber\\
p^{k-1} F 
\left(\begin{smallmatrix} p\tau & pz\\ pz & \tildetau \end{smallmatrix} \right)
+ \frac{1}{p} \sum_{ \lambda \!\!\!\!\!\! \pmod{p}}\!\!\!\!\!
F
\left( \begin{smallmatrix} \frac{ \tau + \lambda}{p} & z\\ z & \tildetau \end{smallmatrix} \right)
\!\!& = &\!\!
p^{k-1} F 
\left(\begin{smallmatrix} \tau & pz\\ pz & p \tildetau \end{smallmatrix} \right)
+ \frac{1}{p} \sum_{ \mu \!\!\!\! \!\!\pmod{p}}  \!\!\!\!\!
F
\left( \begin{smallmatrix} \tau  & z\\ z & \frac{\tildetau + \mu}{p} \end{smallmatrix} \right)\label{drei}\label{3}.
\end{eqnarray}
\end{theorem}
\begin{proof}
We first show that (\ref{eins}) $\Longleftrightarrow$ (\ref{zwei}). 
The direction from left to right is clear since it is a specialization.
The other direction follows from the fact that the Hecke algebra 
$\HH$ is the infinite restricted tensor product of all local Hecke algebras $\HH_p$. 
Here $p$ runs through the set of all primes.
Hence it is sufficient to focus on the generators of $\HH_p$. Here one
has to be careful. 
This conclusion works only because everything is compatible with sums of 
operators and the underlying Hecke algebras are commutative.
Now, since the local Hecke algebras are essentially generated by $T_p$
we are done. 
\\
\\
Next we show that (\ref{zwei}) $\Longleftrightarrow$ (\ref{drei}).
We have seen that $$T_p = \Gamma \left(\begin{smallmatrix} p & 0\\ 0 & 1 
\end{smallmatrix} \right) + \sum_{\lambda \pmod{p}} \Gamma 
\left(\begin{smallmatrix} 1 & \lambda\\ 0 & p \end{smallmatrix} \right).$$
We use this explicit description to calculate $F\vert \widetilde{T_p}^{\bullet}$ and 
$F\vert \widetilde{T_p}_{\bullet}$. Finally we make a change of variable $z \mapsto p^{\frac{1}{2}}z$.
\end{proof}
We parametrize $ Z \in \H_2$ with $Z = \kzxz {\tau}{z}{z}{\tildetau}$. 
Let $X = \kzxz {\tau_x}{z_x}{z_x}{\tildetau_x}$ be the real part of $Z$ and let
$Y= \kzxz {\tau_y}{z_y}{z_y}{\tildetau_y}$ be the imaginary part of $Z$.
Comparing Fourier coefficients in (\ref{3}) we deduce the following result:
\begin{theorem}
Let $k \in \N_0$ be even and let $F:\H_2 \longrightarrow \C$ be a 
$\Gamma$-modular function of weight $k$. Assume that $F$ has
Fourier expansion of the form
\begin{equation}
F(Z) = \sum_{ N } A(N,Y) \, e\{NX\},
\end{equation}
summing over all half-integral symmetric $2\times 2$ matrices.
Then $F\vert_k \bowtie_T = 0$ for all Hecke operators $T \in \HH$ if 
and only if the Fourier coefficients of $F$ satisfy for all prime numbers $p$ the
identity
\begin{equation}
p^{k-1} A
\left( 
\left(\begin{array}{cc}
\frac{n}{p}& \frac{r}{2p} \\ \frac{r}{2p}& m \end{array}\right)
\kzxz {p\tau_y}{pz_y}{pz_y}{\tildetau_y}
\right) + A \left( 
\left(\begin{array}{cc}
pn & \frac{r}{2} \\ \frac{r}{2}& m \end{array}\right)
, 
\kzxz {\frac{\tau_y}{p}}{z_y}{z_y}{\tildetau_y}\right)
\end{equation}
\begin{equation*}
= p^{k-1} A
\left( 
\left(\begin{array}{cc}
n & \frac{r}{2p} \\ \frac{r}{2p}& \frac{m}{p} \end{array}\right)
, \kzxz {\tau_y}{pz_y}{pz_y}{p\tildetau_y}
\right) + A \left( 
\left(\begin{array}{cc}
n& \frac{r}{2} \\ \frac{r}{2}& pm \end{array}\right)
, 
\kzxz {\tau_y}{z_y}{z_y}{\frac{\tildetau_y}{p}}\right).
\end{equation*}
\end{theorem}

\end{document}